\documentclass{amsart}
\usepackage{latexsym}
\usepackage{amsthm,amsfonts,amssymb}
\usepackage[leqno]{amsmath}
\usepackage{enumerate}
\usepackage[all]{xy}\SelectTips{eu}{}
\usepackage{hyperref}

\newcommand{\gr}[1]{\operatorname{gr}_{\fm}(#1)}
\newcommand{\Aplus}[1][0]{\{\rank_RA_i\}_{i\ge #1}}
\newcommand{\Aminus}[1][0]{\{\rank_RA_{i}\}_{i\le #1}}
\newcommand{\kp}{\varkappa}
\newcommand{\Soc}{\operatorname{Soc}}
\newcommand{\lgt}[2][R]{\operatorname{length}_{#1}#2}

\def\ds{\displaystyle}

\def\bA{{\boldsymbol A}}
\def\bP{{\boldsymbol P}}
\newcommand{\BZ}{{\mathbb Z}}

\newcommand{\fm}{{\mathfrak m}}
\newcommand{\fn}{{\mathfrak n}}

\newcommand{\cC}[2]{\operatorname{C}_{#1}{(#2)}}
\newcommand{\PP}{\operatorname{P}}
\newcommand{\II}{\operatorname{I}}
\newcommand{\EE}{\operatorname{E}}

\newcommand{\dd}{\partial}

\newcommand{\HH}[2]{\operatorname{H}_{#1}(#2)}
\newcommand{\CH}[2]{\operatorname{H}^{#1}(#2)}
\newcommand{\rank}{\operatorname{rank}}
\newcommand{\im}{\operatorname{Im}}
\newcommand{\Ker}{\operatorname{Ker}}
\newcommand{\Coker}{\operatorname{Coker}}

\newcommand{\Hom}[3]{\operatorname{Hom}_{#1}({#2},{#3})}
\newcommand{\Ext}[4]{\operatorname{Ext}^{#1}_{#2}(#3,#4)}
\newcommand{\Tor}[4]{\operatorname{Tor}_{#1}^{#2}(#3,#4)}

\newcommand{\xra}{\xrightarrow}
\renewcommand{\le}{\leqslant}
\renewcommand{\leq}{\leqslant}
\renewcommand{\ge}{\geqslant}
\renewcommand{\geq}{\geqslant}

\newcommand{\numberseries}{\bfseries}   

\newlength{\thmtopspace}                
\newlength{\thmbotspace}                
\newlength{\thmheadspace}               
\newlength{\thmindent}                  

\newtheoremstyle{bfupright head,slanted body}
{\thmtopspace}{\thmbotspace}
{\itshape}{\thmindent}{\bfseries}{.}{\thmheadspace}
{{\numberseries \thmnumber{#2.\;}}\thmnote{#3}}

\newtheoremstyle{bfupright head,upright body}
{\thmtopspace}{\thmbotspace}
{\upshape}{\thmindent}{\bfseries}{.}{\thmheadspace}
{{\numberseries \thmnumber{#2.\;}}\thmnote{#3}}

\newtheoremstyle{it head,upright body}
{\thmtopspace}{\thmbotspace}
{\upshape}{\thmindent}{\upshape}{}{\thmheadspace}
{{\numberseries\thmnumber{#2.\!\negthickspace}}
  {}}

\newtheoremstyle{fixed bf head,slanted body}
{\thmtopspace}{\thmbotspace}{\itshape}
{\thmindent}{\bfseries}{.}{\thmheadspace}
{{\numberseries \thmnumber{#2.\;}}\thmname{#1}\thmnote{ (#3)}}

\newtheoremstyle{fixed bf head,upright body}
{\thmtopspace}{\thmbotspace}{\upshape}
{\thmindent}{\bfseries}{.}{\thmheadspace}
{{\numberseries \thmnumber{#2.\;}}\thmname{#1}\thmnote{ (#3)}}

\newtheoremstyle{independent paragraph}
{\thmtopspace}{\thmbotspace}
{\upshape}{\parindent}{\upshape}{}{0pt}
{\thmnote{#3 }}

\newtheoremstyle{subparagraph}
{\thmbotspace}{\thmbotspace}
{\upshape}{\parindent}{\upshape}{}{0pt}
{\thmnote{#3 }}

\theoremstyle{bfupright head,slanted body}
\newtheorem{res}{}[section]             \newtheorem*{res*}{}

\theoremstyle{bfupright head,upright body}
\newtheorem{bfhpg}[res]{}               \newtheorem*{bfhpg*}{}

\theoremstyle{it head,upright body}
\newtheorem{com}[res]{}                 \newtheorem*{com*}{}

\theoremstyle{fixed bf head,slanted body}
          \newtheorem*{thm*}{Theorem}
\newtheorem{prp}[res]{Proposition}      \newtheorem*{prp*}{Proposition}
        \newtheorem*{cor*}{Corollary}
\newtheorem{lem}[res]{Lemma}            \newtheorem*{lem*}{Lemma}

\theoremstyle{fixed bf head,upright body}
\newtheorem{obs}[res]{Observation}      \newtheorem*{obs*}{Observation}
\newtheorem{rmk}[res]{Remark}           \newtheorem*{rmk*}{Remark}

\theoremstyle{independent paragraph}

\theoremstyle{subparagraph}
\newtheorem{spg}{}

\newlength{\thmlistleft}        
\newlength{\thmlistright}       
\newlength{\thmlistpartopsep}   
\newlength{\thmlisttopsep}      
\newlength{\thmlistparsep}      
\newlength{\thmlistitemsep}     

\newenvironment{thmquote}[1][1em]{\begin{list}{}%
    {\setlength{\leftmargin}{#1}\setlength{\rightmargin}{#1}%
      \setlength{\partopsep}{0pt}%
      \setlength{\topsep}{\thmbotspace}%
      \setlength{\parsep}{0pt}%
      \setlength{\itemsep}{0pt}}
      \item[]}
    {\end{list}}%

\newcounter{eqc} 
\newenvironment{eqc}{\begin{list}{\upshape (\textit{\roman{eqc}})}%
    {\usecounter{eqc}%
      \setlength{\leftmargin}{\thmlistleft}%
      \setlength{\labelwidth}{\thmlistleft}%
      \setlength{\rightmargin}{\thmlistright}%
      \setlength{\partopsep}{\thmlistpartopsep}%
      \setlength{\topsep}{\thmlisttopsep}%
      \setlength{\parsep}{\thmlistparsep}%
      \setlength{\itemsep}{\thmlistitemsep}}}%
  {\end{list}}%

\newcounter{prt}
\newenvironment{prt}{\begin{list}{\upshape (\alph{prt})}%
    {\usecounter{prt}%
      \setlength{\leftmargin}{\thmlistleft}%
      \setlength{\labelwidth}{\thmlistleft}%
      \setlength{\rightmargin}{\thmlistright}%
      \setlength{\partopsep}{\thmlistpartopsep}%
      \setlength{\topsep}{\thmlisttopsep}%
      \setlength{\parsep}{\thmlistparsep}%
      \setlength{\itemsep}{\thmlistitemsep}}}%
  {\end{list}}%

\newcommand{\prtlbl}[1]{{\upshape(#1)}}

\newcommand{\prevprt}[1]{{\upshape(\alph{prt}#1)}}

\newenvironment{prf}[1][Proof]{\begin{proof}[\bf #1]}{\end{proof}}

\numberwithin{equation}{res}

\begin{document}

\date{23 November 2006}

\title[Acyclicity over local rings with radical cube zero]{Acyclicity
  over local rings\\with radical cube zero}

\author[L.~W.~Christensen]{Lars Winther Christensen}

\thanks{L.W.C.\ was partly supported by grants from the Danish Natural
  Science Research Council and The Carlsberg Foundation.}

\address{Department of Mathematics, University of Nebraska, Lincoln,
  Nebraska~68588}

\email{winther@math.unl.edu}

\urladdr{http://www.math.unl.edu/{\tiny $\sim$}lchristensen}

\author[O.~Veliche]{Oana~Veliche}

\address{Department of Mathematics, University of Utah, Salt Lake
  City, Utah~84112}

\email{oveliche@math.utah.edu}

\urladdr{http://www.math.utah.edu/{\tiny $\sim$}oveliche}



\begin{abstract}
  This paper studies infinite acyclic complexes of finitely generated
  free modules over a commutative noetherian local ring $(R,\fm)$ with
  $\fm^3=0$.  Conclusive results are obtained on the growth of the
  ranks of the modules in acyclic complexes, and new sufficient
  conditions are given for total acyclicity.  Results are also
  obtained on the structure of rings that are not Gorenstein and admit
  acyclic complexes; part of this structure is exhibited by every ring
  $R$ that admits a non-free finitely generated module $M$ with
  $\Ext{n}{R}{M}{R}=0$ for a few $n>0$.
\end{abstract}

\maketitle

\section*{Introduction}

In this paper, $R$ is a commutative noetherian local ring with maximal
ideal~$\fm$. Throughout, \emph{module} means finitely generated
module.

A chain complex of $R$-modules
\begin{equation*}
  \bA\ \,=\ \, \cdots \longrightarrow A_{i+1} \xra{\partial_{i+1}} A_{i}
  \xra{\;\:\partial_{i}\;\:} A_{i-1} \xra{\partial_{i-1}} A_{i-2}
  \longrightarrow \cdots
\end{equation*}
is \emph{acyclic} if $\HH{}\bA=0$. The focus of this paper is on
complexes of free modules, so we adopt the convention that an
\emph{acyclic complex} consists of free modules. Such a complex $\bA$
is said to be \emph{totally acyclic} if also the dual complex $\bA^*
=\Hom R{\bA}{R}$ is acyclic.

Over a Gorenstein ring, every acyclic complex is totally acyclic
\cite[(4.1.3)]{LWC}. Moreover, a module over such a ring is
maximal Cohen--Macaulay if and only if it is the cokernel of some
differential in an acyclic complex \cite[thm.~(1.4.8) and
(1.4.9)]{LWC}. Thus, acyclic complexes abound over Gorenstein rings.

Over a ring that is not Gorenstein, a non-trivial acyclic complex need
not even exist. Indeed, this is the case for rings that are Golod and
not Gorenstein \cite[exa.~3.5(2)]{LLAAMr02}. Yet, examples of acyclic
complexes over non-Gorenstein rings do exist, and the ones given in
\cite{AGP-97,LWC,OVl01,YYs03} are, in fact, examples of totally
acyclic complexes. It has proved harder to come by acyclic complexes
that are not totally acyclic. However, in \cite{DAJLMS} Jorgensen and
\c{S}ega construct an acyclic, but not totally acyclic, complex over a
local ring with $\fm^3=0$.

\begin{spg}
  This paper started from the observation that the ring considered in
  \cite{DAJLMS} has a specific structure, described by Yoshino
  \cite{YYs03} in a related context. To explain this we introduce some
  notation:
  
  Let $k$ denote the residue field $R/\fm$. Two principal invariants
  of $R$ are the embedding dimension and the socle dimension:
  \begin{equation*}
    e=\rank_k \fm/\fm^2 \qquad\text{and}\qquad r=\rank_k(0:\fm).
  \end{equation*}
  
  The $i$th Bass number of $R$ is $\mu_R^i=\rank_k\Ext iRkR$; note that
  $\mu_R^0 = r$, as $(0:\fm) \cong \Hom{R}{k}{R}$. For an $R$-module
  $M$, the $i$th Betti number is $\beta_i^R(M)=\rank_k\Ext iRMk$. The
  formal power series
  \begin{equation*}
    \II_R(t) = \sum_{i=0}^\infty\mu_R^it^i \qquad\text{and}\qquad 
    \PP_k^R(t) = \sum_{i=0}^\infty\beta_i^R(k)t^i
  \end{equation*}
  are, respectively, the Bass series of $R$ and the Poincar\'{e}
  series of $k$.
  
  A complex of free $R$-modules $\bA$ is \emph{minimal} if
  $\partial(\bA)\subseteq\fm\bA$. In particular, if $\bA$ is minimal
  and acyclic, then either $A_i\ne 0$ for all $i\in\BZ$ or $\bA$ is
  the zero complex. Every acyclic complex contains a minimal one as a
  direct summand with contractible complement.
\end{spg}

\begin{spg}
  In \cite{YYs03} Yoshino proves that when a non-Gorenstein local ring
  $R$ with $\fm^3=0$ admits a non-zero minimal totally acyclic
  complex, either of the two numbers $e$ and $r$ completely determines
  the homological invariants $\II_R(t)$ and $\PP_k^R(t)$. The same
  holds for the ring considered by Jorgensen and \c{S}ega in
  \cite{DAJLMS}.
  
  Let $(R,\fm)$ be a local ring with $\fm^3=0\ne\fm^2$. Suppose $R$ is
  not Gorenstein and admits a non-zero minimal acyclic complex $\bA$.
  This paper considers the following questions:
  \begin{enumerate}[\rm\quad A.]
  \item Does the existence of $\bA$ impose conditions on the structure
    of $R$?
  \item What is the asymptotic behavior of the sequences
    \begin{equation*}
      \Aplus \quad\text{and}\quad \Aminus?
    \end{equation*}
  \item When is $\bA$ totally acyclic?
  \end{enumerate}
  Accordingly, the main results are collected in three theorems.
\end{spg}

\begin{res*}[Theorem A]
  \label{growth}
  Let $(R,\fm,k)$ be a local ring that is not Gorenstein and has
  $\fm^3=0\not=\fm^2$. If there exists a non-zero minimal acyclic
  complex $\bA$ of finitely generated free $R$-modules, then the ring
  has the following properties:
  \begin{prt}
  \item $(0:\fm)=\fm^2.$
  \item $e=r+1$; in particular, $\lgt[]{R} = 2e$.
  \item $\PP _k^R(t) =\ds\frac{1}{(1-t)(1-rt)}$.
  \item[\prevprt{$'$}] The graded ring $\gr{R}$ is Koszul.\footnote{
      $R$ is itself standard graded, i.e.\ $R \cong \gr{R}$, if and
      only if $k$ is a retract of $R$ as a ring, see
      \cite[prop.~1.1]{CLf86}, and by Cohen's Structure Theorem this
      happens if and only if $R$ is equicharacteristic.}
  \end{prt}  
  
  \noindent If, in addition, $\CH{n}{\bA^*}=0$ for some
  integer $n$, then

  \begin{prt}
  \item[\prtlbl{d}] $\II_R(t)=\ds\frac{r-t}{1-rt}$.
  \end{prt}
\end{res*}

\begin{spg}
  Yoshino proved in \cite[thm.~3.1]{YYs03} that $R$ has this structure
  if it is standard graded and $\bA$ is totally acyclic; see also
  Observation~\ref{obs}.
\end{spg}

\begin{spg}
  For the acyclic complex constructed in \cite{DAJLMS}, the sequence
  $\Aplus$ is strictly increasing and $\Aminus$ is constant. A natural
  question, posed in \cite{DAJLMS}, is whether the opposite behavior,
  namely $\Aplus$ constant and $\Aminus$ strictly increasing, is
  possible.  For rings with $\fm^3=0$ the answer is negative:
\end{spg}

\begin{res*}[Theorem B]
  Let $(R,\fm,k)$ be a local ring that is not Gorenstein and has
  $\fm^3=0\not=\fm^2$. If $\bA$ is a non-zero minimal acyclic complex
  of finitely generated free $R$-modules, then one of the following
  holds:
  \begin{enumerate}[\rm(I)]
  \item The residue field $k$ is not a direct summand of $\Coker
    \partial_{i}$ for any $i\in\BZ$, and there is a positive
    integer $a$ such that
    \begin{equation*}
      a = \rank_RA_{i} \quad\text{for all $i\in\BZ$}.
    \end{equation*}
    Moreover, $\operatorname{length}_R \Coker{\dd_i} =ae$ for all
    $i\in\BZ$.
    
  \item There is an integer $\varkappa$, such that $k$ is a direct
    summand of $\Coker \partial_{\varkappa+2}$ and not of $\Coker
    \partial_{i+2}$ for any $i<\varkappa$, and a positive integer $a$
    such that
    \begin{align*}
      a = \rank_RA_{i} \quad &\text{for all integers}\quad i
      \leq\varkappa \quad\text{and}\\
      \rank_RA_{i+1} > \rank_RA_{i} \quad &\text{for all integers}\quad
      i \geq \varkappa.
    \end{align*}
    Moreover, $\operatorname{length}_R \Coker{\dd_{i+2}} =ae$ for all
    integers $i\leq\varkappa$.
  \end{enumerate}
\end{res*}

\begin{spg}
  In case (II) the sequence $\{\rank_R A_i\}_{i\geq \kp}$ has
  exponential growth by work of Lescot \cite[thm.~B]{JLs85}. More
  precise statements about the growth of this sequence are obtained by
  Gasharov and Peeva in \cite[cor.~2.3(ii)]{VNGIVP90} and
  \cite[prop.~3]{IVP98}.
  
  The totally acyclic complex constructed in
  \cite[prop.~3.4]{VNGIVP90} is of type (I), and the acyclic complex
  from \cite[lem.~1.4]{DAJLMS} is of type (II) with  $\kp=0$.
\end{spg}

\begin{spg}
  If $R$ is Gorenstein and $\fm^3=0\ne\fm^2$, then any acyclic complex
  $\bA$ is totally acyclic, and the sequences $\Aplus$ and $\Aminus$
  have the same growth, either exponential or polynomial of the same
  degree. This follows from work of Sj\"{o}din \cite{GSj79}, Lescot
  \cite{JLs85}, and Avramov and Buchweitz \cite{LLAROB00}; see
  \ref{summary} for a summary. If $R$ is not Gorenstein, the first
  implication in Theorem~C contains the result from
  \cite[thm.~3.1]{YYs03} that all modules in a totally acyclic complex
  have the same rank.
\end{spg}

\begin{res*}[Theorem C] Let $(R,\fm)$ be a local ring that is not
  Gorenstein and has $\fm^3=0\not=\fm^2$.  Suppose $\bA$ is a non-zero
  minimal acyclic complex of finitely generated free $R$-modules. Set
  \begin{equation*}
    \mathcal{H}=\{i\in\BZ\mid \CH i{\bA^*}=0\}
  \end{equation*}
  and consider the conditions:
  \begin{eqc}
  \item The set $\mathcal{H}$ contains infinitely many positive
    integers.
  \item All the free modules $A_i$ have the same rank.
  \item If $l-1$ and $l+1$ are in $\mathcal{H}$, then so is $l$.
  \end{eqc}
  The following implications hold:
  \begin{equation*}
    (i)\implies (ii) \implies (iii).
  \end{equation*}
  In particular, if two out of every three consecutive integers belong
  to $\mathcal{H}$, then $\mathcal{H}=\BZ$, i.e.\ $\bA$ is totally
  acyclic.
\end{res*}

\begin{spg}
  This theorem compares to \cite[prop.~2.1]{DAJLMS}, which holds for
  standard graded artinian rings: An acyclic complex $\bA$ is totally
  acyclic if $\BZ \setminus \mathcal{H}$ is a finite set of integers
  of the same parity.
\end{spg}


\section{Modules over local rings with $\fm^3=0$}
\label{modules over rings}

In the rest of this paper, the local ring $(R, \fm,k)$ is assumed to
have $\fm^3=0\ne\fm^2$.  Resolutions of modules over such rings were
first studied by Sj\"{o}din \cite{GSj79} and Lescot \cite{JLs85}; we
open this section with a collection of results from \cite{JLs85}.

In the sequel, the socle $(0:\fm)$ is denoted $\Soc R$. It is clear
that $\fm^2 \subseteq \Soc R$; most of the results from \cite{JLs85}
require $\Soc R=\fm^2$, which is equivalent to assuming $k$ is not a
direct summand of $\fm$; cf.~\cite[lem.~3.2]{JLs85}. This condition is
fulfilled automatically for the rings we are interested in, see
Theorem A.  In fact, it is not too restrictive either: It is not hard
to check that if $k$ is a retract of $R$ as a ring, then $\Soc
R=\fm^2$ or $R$ is a trivial extension of a ring with that property.
That is, $R=Q\ltimes V$, were $(Q,\fn)$ is a local ring with $\fn^3=0$
and $\Soc Q =\fn^2$, and $V$ is a $k$-vector space.

\begin{com}
  \label{socle} From \cite[prop.~3.9(2), thm.~B, and lem.~3.5]{JLs85}
  one has:
  \begin{prt}
  \item If $\Soc R\ne\fm^2$, then for every non-free $R$-module $M$
    the sequence $\{\beta_i^R(M)\}_{i\geq 1}$ is strictly increasing.
  \item If $\Soc R =\fm^2$, then for every $R$-module $M$ the sequence
    $\{\beta_i^R(M)\}_{i\geq 1}$ is eventually constant or has
    exponential growth. In the latter case there is an integer $j$
    such that the sequence $\{\beta_i^R(M)\}_{i\geq j}$ is strictly
    increasing.
  \item If $\Soc R =\fm^2$ and $M$ is an $R$-module with $\fm^2M=0$,
    then
    \begin{equation*}
      \beta_i^R(M) \ge e\beta_{i-1}^R(M) - r\beta_{i-2}^R(M)\ 
      \text{ for all}\ i \ge 2.
    \end{equation*}
  \end{prt} 
\end{com}

\begin{spg}
  The next lemma complements \ref{socle}(b) and contains a special
  case of a result by Gasharov and Peeva \cite[cor.~2.3]{VNGIVP90}.
\end{spg}

\begin{lem}
  \label{betti inequality}
  Assume $\Soc R =\fm^2$ and $e\geq 1+r$. For a non-zero $R$-module
  $M$ with $\fm^2M=0$ there exist integers $ m \ge n \geq 0$, where
  possibly $m=\infty$, such that
  \begin{equation*}
    \cdots > \beta_{m+1}^R(M) > \beta_{m}^R(M) = \cdots =
    \beta_n^R(M) < \beta_{n-1}^R(M) < \cdots < \beta_0^R(M).
  \end{equation*}
  Moreover, if $e>1+r$ then $m=n$ or $m=n+1$.
\end{lem}

\begin{prf} 
  There exists a least $n\geq 0$ such that $\beta_{n+1}^R(M)\geq
  \beta_n^R(M)$. The first inequality below is by \ref{socle}(c),
  \begin{equation*}
    \tag{1}
    \begin{split}
      \beta_{n+2}^R(M) &\ge e\beta_{n+1}^R(M) - r\beta_{n}^R(M)\\
      &\ge e\beta_{n+1}^R(M) - r\beta_{n+1}^R(M)\\
      &= (e-r)\beta_{n+1}^R(M)\\
      &\geq\beta_{n+1}^R(M).
    \end{split}
  \end{equation*}
  By iteration, one has $\beta_{i+1}^R(M)\geq\beta_{i}^R(M)$ for all
  $i\geq n$, and it is immediate that
  $\beta_{j+1}^R(M)>\beta_{j}^R(M)$ implies
  $\beta_{i+1}^R(M)>\beta_{i}^R(M)$ for all $i\ge j$. With
  \begin{equation*}
    m=\inf\{i\in\BZ\mid \beta_{i+1}^R(M) > \beta_{i}^R(M)\} \ge n
  \end{equation*}
  one has $m=\infty$ or $\beta_{i+1}^R(M)>\beta_{i}^R(M)$ for all
  $i\geq m$.
  
  Finally, if $e>r+1$ then $(1)$ yields
  \begin{equation*}
    \beta_{n+2}^R(M) > \beta_{n+1}^R(M),
  \end{equation*}
  which forces $n+1 \ge m \ge n$.
\end{prf}

\begin{spg}
  For an $R$-module $M$, let $M_i$ denote the $i$th syzygy of $M$.
\end{spg}

\begin{com}
  \label{e ge mu} 
  Assume $\Soc R=\fm^2$. Let $M$ be a non-zero $R$-module with
  $\fm^2M=0$, and let $h$ be a positive integer. Following
  \cite[def.~3.1]{JLs85}, $M$ is said to be \emph{$h$-exceptional} if
  $k$ is not a direct summand of the syzygies $M_i$ for $1\leq i \leq
  h$.  If $M$ is $h$-exceptional for every $h\geq 1$, then $M$ is said
  to be \emph{exceptional}.
  
  Let $h\geq 1$. By the proof of \cite[lem.~3.3]{JLs85} an $R$-module
  $M$ is $h$-exceptional if and only if the Betti numbers satisfy:
  \begin{equation}
    \label{eq:e ge mu} 
    \begin{split}
      \beta_1^R(M)&=e\beta_0^R(M)-\rank_k\fm M \quad\text{and}\\
      \beta_i^R(M)&=e\beta_{i-1}^R(M)-r\beta_{i-2}^R(M) \quad\text{for
        all}\quad 2\leq i\leq h.
    \end{split}
  \end{equation}
\end{com}

\begin{com}
  \label{exceptional module} 
  Assume $\Soc R=\fm^2$ and let $h$ be a positive integer. If $R$
  admits an $h$-exceptional module, then $k$ is $h$-exceptional. This
  is \cite[lem.~3.6]{JLs85}.
  
  The equalities \eqref{eq:e ge mu} can be rewritten as an equality of
  polynomials \cite[lem.~3.3]{JLs85}: $k$ is $h$-exceptional if and
  only if
  \begin{equation}
    \label{betti-polynomial}
    \left[\PP_k^R(t)\right]_{\leq
      h}=\left[\ds\frac{1}{1-et+rt^2}\right]_{\leq h},
  \end{equation}
  where $[-]_{\leq h}$ denotes the terms of degree at most $h$.  In
  particular, $k$ is exceptional if and only if
  \begin{equation}
    \label{betti-series}
    \PP_k^R(t)=\ds\frac{1}{1-et+rt^2}.
  \end{equation}
\end{com}

\begin{lem}
  \label{exceptional}
  Assume $R$ is not Gorenstein. If there exists a syzygy module
  $N\ne0$ with $\Ext {h}RNR=0$ for some $h\ge 2$, then $\Soc R =
  \fm^2$.
  
  Moreover, the following hold for an $R$-module $M\ne 0$ with
  $\fm^2M=0$:
  \begin{prt}
  \item If $\,\Ext{2}RMR=0$ and $M$ is a syzygy, then $M$ is
    1-exceptional.
  \item If $\,\Ext{h+1}RMR=0$ for some $h\geq 2$, then $M$ is
    $h$-exceptional.
  \item If $\,\Ext {h}RMR=0$ for infinitely many $h\ge 1$, then $M$ is
    exceptional.
  \end{prt}
\end{lem}

\begin{prf} 
  If $\Soc R\ne\fm^2$, then $k$ is a direct summand of $N_1$ because
  it is a second syzygy; see \cite[lem.~3.2 and proof of
  lem.~3.3]{JLs85}. Therefore,
  \begin{equation*}
    0=\Ext{h}RNR=\Ext {h-1}R{N_1}R =\Ext{h-1}R{k}R\oplus\Ext{h-1}R{N_1'}R
  \end{equation*}
  for some module $N_1'$. In particular, $\Ext{h-1}R{k}R=0$ and that
  contradicts the assumption that $R$ is not Gorenstein.
  
  Note that by this argument, the hypotheses in parts (a)--(c) ensure
  that $\Soc R = \fm^2$, so it makes sense to speak about
  exceptionality.
  
  (a): Applied to $N=M$ and $h=2$ the argument above shows that $M$ is
  $1$-exceptional.
  
  (b): Suppose $M$ is not $h$-exceptional, then $M_i = k\oplus M_i'$
  for some $1\le i \le h$ and some $R$-module $M_i'$. Now there are
  equalities 
  \begin{align*}
    \Ext {h+1-i}RkR \oplus \Ext {h+1-i}R{M_i'}R&=\Ext {h+1-i}R{M_i}R\\
    &=\Ext {h+1}RMR\\
    &=0.
  \end{align*}
  Again, this contradicts the assumption on $R$.  Whence, $M$ is
  $h$-exceptional.
  
  (c): In view of (b), it is sufficient to remark that if $M$ is
  $h$-exceptional for some $h\ge 1$, then $M$ is $i$-exceptional for
  all $1\leq i\leq h$.
\end{prf}

The next result is extracted from the proof of \cite[lem.~3.3]{JLs85}.

\begin{com}
  \label{1-exceptional}
  Assume $\Soc R=\fm^2$.  If $M$ is a $1$-exceptional module, then
  \begin{equation*}
    \rank_k \fm M_1 = r\beta_0^R(M).
  \end{equation*}
\end{com}

\begin{spg}
  For an $R$-module $M$, set $M^* = \Hom RMR$ and let $\ell(M)$ denote
  the length of $M$.  The following equalities are proved already in
  \cite{YYs03}, ostensibly under stronger hypotheses. For
  completeness, we include a proof.
\end{spg}

\pagebreak
\begin{lem} 
  \label{length}
  For an $R$-module $M$ with $\fm^2M=0$ one
  has:
  \begin{prt}
  \item $\ell(M) = \rank_k \fm M + \beta_0^R(M).$
  \item If $\,\Ext{1}R{M}{R}=0$, then $\ell(M^*) =
    r\ell(M)-\beta_0^R(M)\mu_R^1$.
  \end{prt}
\end{lem}

\begin{prf}  
  We can assume $M\ne0$. Set
  \begin{equation*}
    s=\rank_k\fm M \qquad\text{and}\qquad b=\beta_0^R(M).
  \end{equation*}
  
  (a): The exact sequence $ 0 \to \fm M \to M \to M/\fm M \to 0$ is
  isomorphic to
  \begin{equation*}
    \tag{1}
    0 \to k^s \to M \to k^b \to 0.
  \end{equation*}
  In particular, we have $\ell(M) = s + b$.
  
  (b): Dualizing $(1)$, we obtain the exact sequence
  \begin{equation*}
    0 \to {\Hom RkR}^b \to M^* \to {\Hom RkR}^s \to {\Ext 1RkR}^b \to 0,
  \end{equation*}
  which by part (a) and additivity of length yields
  \begin{equation*}
    \ell(M^*) = br + sr - b\mu_R^1 = r\ell(M) - b\mu_R^1.\qedhere
  \end{equation*}
\end{prf}

\begin{lem}
  \label{bass series}
  Assume $R$ is not Gorenstein. If there exists a non-free $R$-module
  $M$ such that $\Ext{n+1}RMR=0$ for some $n\geq 2$, then the Bass
  series of $R$ satisfies:
  \begin{equation*}
    [\II_R(t)]_{\leq n}=\left[\ds\frac{r-et+t^2}{1-et+rt^2}\right]_{\leq
      n}.
  \end{equation*}
\end{lem}

\begin{prf}
  It follows from Lemma~\ref{exceptional} that $\Soc R = \fm^2$.  Let
  $E=\EE(k)$ be the injective envelope of $k$.  Since the module $M$
  is not free and
  \begin{equation*}
    \Tor{n+1-i}{R}{M}{E_i} = \Tor{n+1}RME=\Hom R{\Ext{n+1}RMR}{E}=0  
  \end{equation*}
  for $0 < i \leq n$, the syzygy module $E_1$ is $n-1$ exceptional
  and does not contain $k$ as a direct summand.  Now
  \cite[2.8(3\&4)]{HSV-04} yield
  \begin{equation*}
    \tag{1}
    \beta^R_0(E_1) = e(r-1) \quad\text{and}\quad \ell(E_1)=(r-1)(1+e+r),
  \end{equation*}
  and from Lemma~\ref{length}(a) we get
  \begin{equation*}
    \tag{2}
    \rank_k \fm E_1 = (r-1)(1+e+r)-e(r-1)=r^2-1.
  \end{equation*}
  The Betti numbers of the module $E$ are the Bass numbers of $R$;
  that is $\mu_R^i = \beta^R_i(E) = \beta^R_{i-1}(E_1)$. Rewriting the
  equations \eqref{eq:e ge mu} for the module $E_1$ as an equality of
  polynomials gives the first equality below. The second follows by
  $(1)$ and $(2)$.
  \begin{align*}
    [\II_R(t)]_{\leq n} &= r + t\left[\ds\frac{\beta^R_0(E_1) -
        (\rank_k\fm
        E_1)t}{1-et+rt^2}\right]_{\leq n-1}\\
    &=r+t\left[\ds\frac{e(r-1)- (r^2-1)t}{1-et+rt^2}\right]_{\leq n-1}\\
    &=\left[\ds\frac{r(1-et+rt^2)+e(r-1)t-
        (r^2-1)t^2}{1-et+rt^2}\right]_{\leq n}\\
    &=\left[\ds\frac{r-et+t^2}{1-et+rt^2}\right]_{\leq n}.\qedhere
  \end{align*} 
\end{prf}


\section{Proofs of Theorems A-C}
\label{Rings with acyclic complexes}

In this section we prove the three main theorems, stated in the
Introduction.

\begin{com}
  \label{notation}
  Let $\bA$ be a minimal acyclic complex; throughout this section we
  use the following notation:
  \begin{equation*}
    b_i = \rank_R A_i,\ \
    \cC{i}\bA=\Coker\dd_{i+1} \cong \Ker\dd_{i-1}, \ \text{ and } \ s_i =
    \rank_k{\fm\!\cC{i}\bA}
  \end{equation*}
  for $i\in\BZ$. Note that $\beta_0^R(\cC{i}\bA) = b_i$.
\end{com}

\begin{rmk}
  \label{k summand}
  Assume $R$ is not a hypersurface ring, the Betti numbers of $k$ are
  then strictly increasing; see \cite[rem.~8.1.1(3)]{LLA98}.  Let
  $\bA$ be an acyclic $R$-complex; note that for any $j\in\BZ$ the
  inequalities
  \begin{equation*}
    \beta^R_{b_j}(k) > \beta^R_{b_j-1}(k) > \cdots >
    \beta^R_{0}(k) = 1
  \end{equation*}
  show that $\beta^R_{b_j}(k) > b_j$, so $k$ cannot be a direct
  summand of $\cC{i}{\bA}$ for any $i\le j-b_j$.  In particular, $k$
  is not a direct summand of $\cC{i}{\bA}$ for any $i\le -b_0$.
  
  For a minimal acyclic complex $\bA$ set
  \begin{equation*}
    \kp =\inf\{i\mid k\ \text{\rm is a direct summand of } \cC{i+1}{\bA} \}
  \end{equation*}
  and note that
  \begin{equation*}
    \infty\geq \kp \geq -b_0 >-\infty.
  \end{equation*}
  If $\Soc R=\fm^2$, then \eqref{eq:e ge mu} yields
  \begin{align}
    \label{eq betti}
    b_i &= eb_{i-1} -rb_{i-2} \quad\text{for all}\quad i \le \kp,\quad
    \text{and}\\
    \label{ieq betti}
    b_{\kp+1} &> eb_{\kp} -rb_{\kp-1}.
  \end{align}
\end{rmk}

\begin{bfhpg}[Proof of Theorem~A] We may, after a shift, 
  assume $k$ is not a direct summand of $\cC{-i}{\bA}$ for any $i\ge
  0$; cf.~Remark~\ref{k summand}.
  
  (a): Suppose $\Soc R\not=\fm^2$; by \ref{socle}(a) one gets
  \begin{equation*}
    b_0 > b_{-1} > b_{-2}>\cdots>b_{-b_0}>0,
  \end{equation*}
  which is absurd. Therefore, $\Soc R=\fm^2$.
  
  (b): Set $a_i=b_{-i}$, then \eqref{eq betti} translates to
  \begin{equation*}
    a_{i} = ea_{i+1} - ra_{i+2} \quad\text{for all}\quad i\geq 0.
  \end{equation*}
  By Proposition~\ref{infinite sequence}(b) it follows that $e=r+1$.
  
  (c): For $h\ge 1$ our assumption on the complex $\bA$ implies that
  $\cC{-h}\bA$ is an $h$-exceptional $R$-module. By \ref{exceptional
    module}, $k$ is then $h$-exceptional for all $h\geq 1$ and hence
  exceptional. The expression for the Poincar\'{e} series of $k$ now
  follows from \eqref{betti-series} and part (b).
  
  (c$'$): It follows from (b) and (c) that the Hilbert series of
  $\gr{R}$ is
  \begin{equation*}
    1+et+rt^2 = (1+t)(1+rt) = \frac{1}{\PP_k^R(-t)}.
  \end{equation*}
  As $\PP_k^R(t) = \PP_k^{\gr{R}}(t)$ by \cite[thm.~2.3]{CLf86}, it
  follows that $\gr{R}$ is a Koszul algebra by \cite[thm.~1.2]{CLf86}.
  
  (d): After a shift, we may assume $n =-1$; then $\CH{n}{\bA^*}=0$
  translates to
  \begin{equation*}
    \Ext {i+1}R{\cC{-i}\bA}R=0 
    \quad\text{for all}\quad i\geq 0.
  \end{equation*}
  Lemma~\ref{bass series} now applies to the modules $\cC{-i}\bA$ for
  all $i > 0$ and yields
  \begin{equation*}
    \II_R(t) = \frac{r-et+t^2}{1-et+rt^2}.
  \end{equation*}
  Using part (b) we obtain, after simplification, the desired
  equality.\qed
\end{bfhpg}

\begin{bfhpg}[Proof of Theorem~B]  
  \pushQED{\qed}%
  Note that $\Soc R = \fm^2$ and $e=r+1$ by Theorem~A. It was already
  remarked in \ref{k summand} that $\infty\geq \kp > -\infty$. Thus,
  either $\kp=\infty$ or we may, after a shift, assume $\kp=0$.  The
  first case corresponds to (I) and the second to (II).  Set
  $a_i=b_{-i}$; in either case \eqref{eq betti} translates to
  \begin{equation*}
    a_{i} = ea_{i+1} - ra_{i+2} \quad\text{for all}\quad i\geq 0.
  \end{equation*}
  By Proposition~\ref{infinite sequence}(a) there is a positive
  integer $a$, such that
  \begin{equation*}
    b_{i} =a \quad\text{for all}\quad i\le 0.
  \end{equation*}
  
  In case $\kp =0$, the inequality \eqref{ieq betti} becomes
  \begin{equation*}
    b_{1}-b_0 > r(b_0-b_{-1})=0.
  \end{equation*}
  Thus, $b_{1}>b_0$ and Lemma~\ref{betti inequality} applied to
  $\cC{0}\bA$ yields the desired conclusion.
  
  In case $\kp = \infty$, the equality \eqref{eq betti} translates to
  \begin{equation*}
    b_{i}-b_{i-1} = r(b_{i-1}-b_{i-2}) \quad\text{for all}\quad i\in\BZ.
  \end{equation*}
  Since $b_0 = a = b_{-1}$, it follows by recursion that $b_i=a$ also
  for $i > 0$.
  
  For $i\le\kp$, the residue field $k$ is not a direct summand of
  $\cC{i}{\bA}$. Therefore, one has $\rank_k\fm\!\cC{i}{\bA} =ar$ by
  \ref{1-exceptional}, and Lemma~\ref{length}(a) yields the desired
  \begin{equation*}
    \ell(\cC{i}{\bA}) = ar + a = ae.\qedhere
  \end{equation*}
\end{bfhpg}

\begin{bfhpg}[Proof of Theorem~C]
  \pushQED{\qed}%
  First note that $\Soc R = \fm^2$ and $e=r+1$ by Theorem~A.
  
  $(i)\!\!\implies\!\!(ii)$: Let $C$ be any cokernel $\cC{i}{\bA}$. By
  assumption, $\Ext{h}RCR = 0$ for infinitely many $h>0$, so $C$ is
  exceptional by Lemma~\ref{exceptional}(c).  Thus, $k$ is not a
  direct summand of any cokernel $\cC{i}{\bA}$, and it follows by
  Theorem~B that all the modules $A_i$ have the same rank.
  
  $(ii)\!\!\implies\!\!(iii)$: After a shift we may assume $l=1$, so
  $0$ and $2$ are in $\mathcal{H}$.  Consider the dual complex
  \begin{equation*}
    \bA^* = \dots \to A^*_{-2} \xra{\dd_{-1}^*} A^*_{-1} \xra{\dd_0^*}
    A^*_0 \xra{\dd_1^*} A^*_1 \xra{\dd_2^*} A^*_2 \to \cdots.
  \end{equation*}
  By definition, $\CH {1}{\bA^*} = \Ker \dd_2^* / \im \dd_1^*$.  We
  will show that $\Ker \dd_2^*$ and $\im \dd_1^*$ have the same
  length, $ae$, where $a$ is the common rank of the modules $A_i$.
  
  First note that for all $i\in\BZ$ we have $\Ker \dd_i^*=(\Coker
  \dd_i)^*={\cC{i-1}\bA}^*$, by left-exactness of $\Hom{R}{-}{R}$.  By
  assumption $\CH 0{\bA^*} = 0 = \CH {2}{\bA^*}$; this means that
  $\Ext 1R{\cC{-1}\bA}R=0$ and $\Ext 1R{\cC{1}\bA}R=0$. For $i=\pm 1$
  it follows by Lemma \ref{length}(b), Theorem~B, and Theorem~A(b,d)
  that
  \begin{align*}
    \ell({\cC{i}\bA}^*) &= r\ell(\cC{i}\bA) - \beta_{0}^R(\cC{i}\bA)\mu_R^1\\
    &= ra(r+1) - a(r^2-1)\\
    &= a(r+1) \\
    &= ae.
  \end{align*}  
  Thus we have
  \begin{equation*}
    \ell(\Ker \dd_2^*)=ae = \ell(\Ker \dd_0^*).
  \end{equation*}
  For all $i\in\BZ$ we have $\ell(\Ker \dd_i^*) + \ell(\im \dd_i^*) =
  a\ell(R)$; moreover, since $\CH 0{\bA^*} = 0$ we have $\ell(\Ker
  \dd_1^*) = \ell(\im \dd_0^*)$. Combining these equations we find
  \begin{align*}
    \ell(\im \dd_1^*) &= a\ell(R) - \ell(\Ker \dd_1^*)\\
    &= a\ell(R) - \ell(\im \dd_0^*)\\
    &= \ell(\Ker \dd_0^*)\\
    &= ae. \qedhere
  \end{align*}
\end{bfhpg}

\section{Concluding remarks and questions}

In this section we sum up the state of the three questions raised in
the Introduction. The assumption that $(R,\fm,k)$ is local with
$\fm^3=0\ne\fm^2$ is still in force.

\subsection*{A. Structure of a non-Gorenstein ring  admitting an
  acyclic complex} One answer to this question is given by Theorem~A.
It remains open whether the additional assumption, in Theorem~A(d),
that some cohomology module vanishes, is fulfilled automatically. See
also Question~\ref{Q1} below.

\begin{spg}
  It also remains open whether every non-Gorenstein ring $R$ with the
  structure described in Theorem~A admits a non-zero minimal acyclic
  complex. For a construction of totally acyclic complexes over
  certain rings, see \cite[thm.~(3.1)]{AGP-97}.
  
  If one allows for non-finitely generated modules, an acyclic
  $R$-complex can always be constructed by copying part of the
  argument for \cite[prop.~6.1(3)]{SInHKr}: Let $\bP$ be a projective
  resolution of the injective hull of $k$, then the mapping cone of
  the homothety morphism $R \to \Hom{R}{\bP}{\bP}$ is an acyclic
  complex of flat $R$-modules, and flat modules are free, as $R$ is
  artinian.
\end{spg}

\subsection*{B. Asymptotic behavior of ranks}
For a minimal acyclic $R$-complex $\bA$, the possible asymptotic
behaviors of $\Aplus$ and $\Aminus$ are now completely understood. For
non-Gorenstein rings it is explained by Theorem~B. For Gorenstein
rings we collect the results in:

\begin{bfhpg}[Summary]
  \label{summary}
  Let $R$ be Gorenstein, then $\Soc R = \fm^2$.  For a minimal acyclic
  complex $\bA$, the sequences $\Aplus$ and $\Aminus$ have the same
  growth, either exponential or polynomial of the same degree. We show
  this below by arguing on the embedding dimension of $R$, but first
  we make an observation: Because $R$ is artinian and Gorenstein,
  every non-free $R$-module is a cokernel in a minimal acyclic
  complex, which is determined uniquely up to isomorphism. By the
  Krull--Schmidt theorem every $R$-module decomposes uniquely as a sum
  of indecomposable modules, and it follows that every minimal acyclic
  complex is isomorphic to a sum of acyclic complexes whose cokernels
  are indecomposable.
  
  \begin{spg}
    If $e=1$, then all modules in a minimal acyclic $R$-complex have
    the same rank. Indeed, by Cohen's Structure Theorem, $R$ is
    isomorphic to $D/t^3D$, where $D$ is a discrete valuation domain
    with maximal ideal $tD$. Up to isomorphism there are, therefore,
    three indecomposable $R$-modules: $R$, $R/\fm$, and $R/\fm^2$. The
    two non-free ones have constant Betti numbers equal to $1$.
  \end{spg}
  
  \begin{spg}
    If $e=2$, then $R$ is a complete intersection ring.  Let $C$ be a
    cokernel in $\bA$. In view of the isomorphism $\Ext{i}{R}{C}{C}
    \cong \Ext{i}{R}{C^*}{C^*}$ it follows by
    \cite[cor.~5.7]{LLAROB00} that $\Aplus$ and $\Aminus$ have
    polynomial growth of the same degree. By Lemma~\ref{betti
      inequality} it follows that either all the modules in $\bA$ have
    the same rank or, after a shift, one has
    \begin{equation*}
      \cdots > \rank_R A_{m+1} > \rank_R A_m = \dots = \rank_R A_0 <
      \rank_R A_{-1} < \cdots,
    \end{equation*}
    where $\infty > m \ge 0$.
    
    The indecomposable modules over the ring $R=k[X,Y]/(X^2,Y^2)$ are
    classified in \cite{AHlIRn61}.  For every $l>0$ an indecomposable
    module of even length $2l$ determines an acyclic complex $\bA$
    with $\rank_R A_n = l$ for all $n\in\BZ$; this follows from
    \cite[prop.~5]{AHlIRn61} and Lemma~\ref{length}(a).  The
    indecomposable modules of odd length are exactly the syzygies and
    cosyzygies of $k$; see also \cite[4.2.3]{LLAROB00a}.  After a
    shift, these modules all determine the same acyclic complex $\bA$,
    for which $\{\rank_RA_i\}_{i \in \BZ}$ is the sequence
    \begin{align*}
      \cdots > n+1 > n > \cdots > 2 > 1 = 1 < 2 < \cdots < n < n+1 <
      \cdots.
    \end{align*}
  \end{spg}
  
  \begin{spg}
    If $e\ge 3$ and $\bA$ is non-zero, then $\Aplus$ and $\Aminus$
    have exponential growth by \cite[thm.~B]{JLs85}. Moreover, after a
    shift one has
    \begin{equation*}
      \cdots > \rank_R A_2 > \rank_R A_1 \ge \rank_R A_0 < \rank_R
      A_{-1} < \cdots 
    \end{equation*}
    by Lemma~\ref{betti inequality}. An example is $R=k[X,Y,Z]/(X^2
    -Y^2, Y^2 - Z^2, XY, YZ)$.
  \end{spg}
\end{bfhpg}

\subsection*{C. Acyclicity and total acyclicity}
All acyclic complexes known to the authors, including the one from
\cite{DAJLMS}, can be obtained by a standard technique:

\begin{bfhpg}[Construction]
  \label{construction}
  Let $M$ be an $R$-module; take minimal free resolutions
  \begin{equation*}
    \xymatrix@R=1ex@C=1.5em{
      \cdots \ar[r] & P_2 \ar[r]^-{d_2} & P_1 \ar[r]^-{d_1} &
      P_0 \ar[r]^-{p} & M \ar[r] & 0\\
      \cdots \ar[r] & Q_2 \ar[r]^-{\dd_2} & Q_1 \ar[r]^-{\dd_1} &
      Q_0 \ar[r]^-{\pi} & M^* \ar[r] & 0}
  \end{equation*}
  and form the complex
  \begin{equation}
    \label{eq:construction}
    \bA=\cdots \to Q_2 \xra{\dd_2} Q_1 \xra{\dd_1} Q_0
    \xra{p^*\circ\pi} P_0^* \xra{d_1^*} P_1^* \xra{d_2^*}
    P_2^* \to \cdots
  \end{equation}
  with $P^*_0$ in degree $0$.  If $\Ext iRMR=0$ for all $i>0$, then
  $\bA$ is acyclic and
  \begin{equation*}
    M = \Coker\dd^{\bA^*}_1 \qquad\text{and}\qquad \CH
    i{\bA^*}=0\quad\text{for all}\quad i<0.
  \end{equation*}
  Moreover, if $M^*$ is without non-zero free direct summands, then
  $\bA$ is minimal.
  
  On the other hand, if $\bA$ is some acyclic complex with $\CH
  i{\bA^*}=0$ for all $i<0$, then the module $M = \Coker\dd^{\bA^*}_1$
  has $\Ext iRMR=0$ for all $i>0$.
\end{bfhpg}

\begin{obs} 
  \label{obs}
  In combination, \cite[prop.~2.9]{HSV-04}\footnote{ Which contains a
    typo: the equalities of Betti numbers should be
    $b_0(M)=b_1(M)=\dots=b_{j+1}(M)$.} and Lemma~\ref{bass series} can
  be reformulated as follows:
  \begin{thmquote}
    \setlength{\leftmargin}{1em} \it Let $(R,\fm, k)$ be a local ring
    that is not Gorenstein and has $\fm^3=0\not=\fm^2$. If there
    exists an $R$-module $M\ne 0$ and an integer $n\geq 3$ such that
    $\fm^2M=0$ and
    \begin{equation*}
      \Ext {n-1}RMR=\Ext{n}RMR=\Ext {n+1}RMR,
    \end{equation*}
    then there are equalities $\beta_{n}^R(M) = \dots = \beta_1^R(M)=
    \beta_0^R(M)$, and the ring has properties \prtlbl{a} and
    \prtlbl{b} from Theorem~A. Moreover, the ring satisfies:
    \begin{prt}
    \item[\prtlbl{c}] $\left[\PP _k^R(t)\right]_{\leq
        n}=\left[\ds\frac{1}{(1-t)(1-rt)}\right]_{\leq {n}}$.
    \item[\prtlbl{d}] $\left[\II_R(t)\right]_{\leq {n}} =
      \left[\ds\frac{r-t}{1-rt}\right]_{\leq {n}}$.
    \end{prt}
  \end{thmquote}
  
  If $M$ is a module with $\Ext iRMR=0$ for all $i>0$, and $\bA$ is
  the corresponding acyclic complex, cf.~\ref{construction}, then
  $\bA$ satisfies the hypothesis of Theorem~A. The module $M$
  satisfies the hypothesis of the result above for all integers $n\geq
  3$, so it yields the same conclusion as Theorem~A.
\end{obs}

This naturally raises the following

\begin{bfhpg}[Question]
  \label{Q1}
  Is every minimal acyclic $R$-complex $\bA$ obtainable from a module
  $M$ with $\Ext iRMR=0$ for all $i>0$ and thus of the form
  \eqref{eq:construction}?
\end{bfhpg}

The authors are not aware of any example of an acyclic, but not
totally acyclic, complex in which all the modules have the same rank.
Hence one may even ask

\begin{bfhpg}[Question]
  For an acyclic $R$-complex $\bA$ set
  \begin{equation*}
    \kp =\inf\{i\mid k\ \text{\rm is  a direct summand of }
    \cC{i+1}{\bA} \}
  \end{equation*}
  as in Remark~\ref{k summand}. Does one always have $\CH i{\bA^*}=0$
  for all $i < \kp$?
\end{bfhpg}

\appendix
\section*{Appendix}
\stepcounter{section}

Here we prove a technical result on sequences of positive integers
that satisfy a certain second order linear recursion formula.
\begin{prp}
  \label{infinite sequence}
  Let $e>0$ and $r>1$ be integers.  If there exists a sequence of
  positive integers $\{a_i\}_{i\geq 0}$ such that
  \begin{equation*}  
    a_i=ea_{i+1}-ra_{i+2}\quad \text{for all}\quad  i\geq 0,
  \end{equation*}
  then the following hold:
  \begin{prt}
  \item The sequence $\{a_i\}_{i\geq 0}$ is constant.
  \item $e=r+1$.
  \end{prt}
\end{prp}

\begin{prf}  
  Set $q_i={a_i}/{a_{i+1}}$. From the recursion formula one gets for
  each $i\ge 0$:
  \begin{equation*}
    \tag*{\ensuremath{(1)_i}} q_i = e- \frac{r}{q_{i+1}}.
  \end{equation*}
  Subtract $(1)_{i+1}$ from $(1)_i$ to get
  \begin{equation*}
    \tag*{\ensuremath{(2)_i}} q_i - q_{i+1} = r\frac{q_{i+1} -
      q_{i+2}}{q_{i+1}q_{i+2}}.
  \end{equation*}
  
  Next we show that
  \begin{equation*}
    \tag*{\ensuremath{(3)}}
    q_i=q_0\quad\text{for all}\quad i\geq 0.
  \end{equation*}
  Let $i\ge 0$; multiplying the equations $(2)_0,\dots,(2)_{i}$ one
  gets
  \begin{equation*}
    q_0-q_1=r^{i+1}\frac{q_{i+1}-q_{i+2}}{q_1q_2^2\cdot\dots\cdot
      q_{i+1}^2 q_{i+2}}.
  \end{equation*}
  Rewrite this equality in terms of the $a_i$s and simplify as follows
  \begin{align*}
    \frac{a_0}{a_1}-\frac{a_1}{a_2} &=r^{i+1} \frac{
      \ds\frac{a_{i+1}}{a_{i+2}} -\frac{a_{i+2}}{a_{i+3}} } {\ds
      \frac{a_1}{a_2} \cdot \frac{a_2^2}{a_3^2} \cdot\dots\cdot
      \frac{a^2_{i+1}}{a^2_{i+2}} \cdot \frac{a_{i+2}}{a_{i+3}} }
    \\[0.5ex]
    &=r^{i+1} \frac{\ds
      \frac{a_{i+1}a_{i+3}-a_{i+2}^2}{a_{i+2}a_{i+3}} } {\ds
      \frac{a_1a_2}{a_{i+2}a_{i+3}} }
    \\[0.5ex]
    &=r^{i+1}\frac{a_{i+1}a_{i+3}-a_{i+2}^2}{a_1a_2}.
  \end{align*}
  Multiplication by $a_1a_2$ yields
  \begin{equation*}
    \tag{4}
    a_0a_2-a_1^2=r^{i+1}(a_{i+1}a_{i+3}-a_{i+2}^2)\quad\text{for all}\quad
    i\geq 0.
  \end{equation*}   
  Thus, $r^{i+1}$ divides $a_0a_2-a_1^2$ for all $i\geq 0$, which
  forces $a_0a_2-a_1^2=0$ as $r>1$. By $(4)$ we now have
  $a_ia_{i+2}-a_{i+1}^2=0$ for all $i\geq 0$; that is, $q_i = q_{i+1}$
  and therefore $q_i=q_0$ for all $i\geq 0$.
  
  The recursion formula may now be rewritten
  \begin{equation*}
    \tag{5}
    q_0^2 - eq_0 + r = 0.
  \end{equation*}
  Since $q_0$ is rational, it follows by the Rational Root Test that
  $q_0$ is an integer and divides $r$. If $q_0>1$, then
  \begin{equation*}
    a_0>a_1>a_2>\cdots >0,
  \end{equation*}
  which is impossible. Thus, $q_0=1$ and then $(3)$ implies part (a)
  while (b) follows from $(5)$.
\end{prf}


\section*{Acknowledgments}

The authors thank Luchezar Avramov and R\u{a}zvan Veliche for
interesting discussions related to this material. Thanks are also due
to Srikanth Iyengar, Greg Piepmeyer, and Diana White for useful
comments on the exposition.


\providecommand{\bysame}{\leavevmode\hbox to3em{\hrulefill}\thinspace}
\renewcommand{\MR}{\relax\ifhmode\unskip\space\fi MR }

\end{document}